\documentstyle[12pt]{article}
\newcommand{\be}{\begin{equation}}
      \newcommand{\ee}{\end{equation}}
      \newcommand{\ba}{\begin{eqnarray}}
       \newcommand{\ea}{\end{eqnarray}}
\newcommand{\ban}{\begin{eqnarray*}}
       \newcommand{\ean}{\end{eqnarray*}}

\renewcommand{\pt}{\partial}

\newcommand{\ra}{\rightarrow}

 \renewcommand{\o}[2]{\frac{#1}{#2}}
\newcommand{\hf}{\o{1}{2}}

 \newcommand{\qed}{\hspace*{\fill}\rule{3mm}{3mm}\quad}
 \newcommand{\Pf}{\noindent {\em Proof.} }

\newcommand{\heta}{\hat{\eta}}
\newcommand{\Tr}{{\rm Tr}}
\newcommand{\ind}{{\rm ind}}
\newcommand{\tB}{\tilde{B}}
\newcommand{\tBt}{\tilde{B}_t}
\newcommand{\tA}{\tilde{A}}
\newcommand{\tP}{\tilde{P}}
\newcommand{\tD}{\tilde{D}}
\newcommand{\End}{{\rm End}}

\newcommand{\spec}{{\rm spec}}
\newcommand{\Res}{{\rm Res}}

\font\BBb=msbm10 at 12pt
\newcommand{\Bbb}[1]{\mbox{\BBb #1}}

\newcommand{\sect}[1]{\section{#1} \setcounter{equation}{0}}

\newtheorem{theo}{Theorem}[section]

\begin{document}
\newtheorem{lem}[theo]{Lemma}
\newtheorem{prop}[theo]{Proposition}  
\newtheorem{coro}[theo]{Corollary}

\title{Higher Spectral Flow}
\author{Xianzhe Dai\thanks{Partially supported by the Alfred  
P. Sloan Foundation} \\
  \and Weiping Zhang\thanks{Partially supported by the Chinese National Science 
Foundation and the State Education Commission of China.
}}
\maketitle

\begin{abstract}
For a continuous curve of families of Dirac type operators we define a higher spectral flow 
as a $K$-group element. We show that this higher spectral flow can be computed 
analytically by $\heta$-forms, and is related to the family index in the same way as the 
spectral flow is related to the index.
We  introduce a notion of Toeplitz family and relate its index to the higher spectral 
flow. Applications to family indices  for manifolds with boundary are also given.
\end{abstract}

{\bf Introduction}

$\ $

The spectral flow for a one parameter family of self adjoint Fredholm operators is an 
integer that counts the {\em net number} of eigenvalues that change sign. This notion is 
introduced by Atiyah-Patodi-Singer \cite{aps} in their study of index theory on manifolds 
with boundary and is intimately related to the $\eta$ invariant, 
which is also introduced by Atiyah-Patodi-Singer [APS2]. It has since then found other
significant applications.

In this paper, a notion of {\em higher spectral flow} is introduced, generalizing  the usual 
spectral flow. The higher spectral flow is defined for a continuous one parameter family 
(i.e. a curve) of {\em families of Dirac type operators} parametrized by a compact space 
and is an element of the $K$-group of the parameter space. The (virtual) dimension of this 
$K$-group element is precisely the (usual) spectral flow. The definition makes use of the 
concept of {\em spectral section} introduced recently by Melrose-Piazza
\cite{mp}.

Roughly speaking, a spectral section is a  way of dividing the spectral data into the positive 
and negative parts. And two different spectral sections  give rise to a difference element 
which lies in the $K$-group of the parameter space. In order to define the spectral flow for 
a curve of families of operators, one first needs to fix a way of dividing the spectral data at the two 
end points of the curve, i.e., two spectral sections. One then continuously deform the 
begining spectral section along the curve. The diffence (element) of the ending spectral 
section with the (other) given spectral section is the (higher) spectral flow.

More precisely, let $\pi:\,  X \longrightarrow B$ be a smooth fibration with the typical fiber 
$Z$ an odd 
dimensional closed manifold and $B$ compact. A family of self adjoint elliptic 
pseudodifferential operators on $Z$, parametrized by $B$, will be called a $B$-family. 
Consider
a curve of $B$-families, $D_u=\{D_{b,u}\}$, $u\in [0,1]$. 

Assume that the index bundle of $D_0$ vanishes and let $Q_0$, $Q_1$ be spectral 
sections of $D_0$, $D_1$ respectively.  If we consider the total family 
$\tilde{D}=\{D_{b, 
u}\}$ parametrized by $B\times I$, then there is also a total spectral section 
$\tilde{P}=\{P_{b,u}\}$. Let $P_u$ be the restriction of $\tilde{P}$ over $B\times 
\{u\}$. The (higher) spectral flow ${\rm sf}\{(D_0, Q_0),\ (D_1, Q_1)\}$ between the pairs 
$(D_0, 
Q_0)$, $(D_1, Q_1)$ is an element in $K(B)$ defined by
\be {\rm sf}\{(D_0, Q_0),\ (D_1, Q_1)\} = [Q_1-P_1] -  [Q_0-P_0] \in K(B). \ee
The definition is independent of the choice of the (total) spectral section $\tilde{P}$.
When $D_u$, $u\in S^1$ is a periodic family, we choose $Q_1=Q_0$. In this case the 
higher
spectral flow turns out to be independent of $Q_0=Q_1$ and therefore defines an invariant 
of the family, denoted by ${\rm sf}\{D_u\}$.

We show that this is indeed a generalization of the notion of spectral flow (Theorem 1.3). 
That is, when we consider a one parameter family of self adjoint elliptic pseododifferential 
operators and the standard Atiyah-Patodi-Singer  spectral projections [APS2]
at the end points, the higher spectral flow coincides with the usual spectral flow. 

We also show that this higher version of spectral flow satisfies the basic properties
of spectral flow. For example, its Chern character can be expressed analytically in terms of 
a generalized version
of the $\heta$ form of Bismut-Cheeger \cite{bc1}. This gives a way of 
computing, analytically, the Chern character of the higher spectral flow.

\begin{theo}
Let $D_u$ be a curve of $B$-families of Dirac type operators and $Q_0$, $Q_1$ spectral 
sections of $D_0$, $D_1$ respectively. Let $\tB_t(u)$ be a curve of superconnections. 
Then
we have the following identity in $H^*(B)$,
\be  
{\rm ch}({\rm sf}\{(D_0, Q_0),\ (D_1, Q_1)\}) = \heta(D_1, Q_1) -\heta(D_0, Q_0)
 -\o{1}{\sqrt{\pi}}  \int^1_0 a_0(u) du , \label{hsf}
\ee
where $\heta(D_i, Q_i)$ are suitable $\heta$-forms  
and  $a_0(u)$ is a local invariant determined by the asymptotic expansion  (2.22).
\end{theo}

In the case of periodic family, the higher spectral flow can be  related to the family index, in 
the same way as the spectral flow is related to the index. Namely,  let $D=\{ D_u\} _{u\in 
S^1}$ be a periodic family of $B$-families of
(self-adjoint) Dirac type operators $\{ D_{b,u}\} _{b\in B,u\in S^1}$.
Then for any $b\in B$, $\{ D_{b,u}\}_{u\in S^1}$ determines a natural Dirac type
operator $D_b^\prime$ on $S^1\times Z_b$:
$$D_b^\prime:\Gamma ({\Bbb C}^2\otimes S(TZ_b)\otimes E_b)\rightarrow \Gamma
({\Bbb C}^2\otimes S(TZ_b)\otimes E_b). $$
$D^\prime =\{ D_b^\prime \}_{b\in B}$ then forms a family of Dirac type operators
over $B$.

\begin{theo} The following identity holds in $H^*(B)$,
$${\rm ch}({\rm ind}\, D^\prime)={\rm ch}({\rm sf}\{ D\}). $$
\end{theo}

This theorem generalizes the well-known result of Atiyah-Patodi-Singer \cite{aps} to the 
family case.  As a consequence we obtain the following multiplicativity formula for the 
spectral flow. 

\begin{coro} Assume that $B$ is also spin and for each $u\in S^1$, denote $D^X_u$ the 
total Dirac operator on $X$. Then the following identity holds,
$${\rm sf}\{ D_{u\in S^1}^X\}=\int_B\hat{A}(TB){\rm ch}({\rm sf}\{ D_{u\in S^1}^Z \}).$$
\end{coro}

 Using spectral sections, we also introduce a notion of Toeplitz family, extending the usual concept of Toeplitz operator to the fibration case.
Namely let $g:X\rightarrow GL(N,{\Bbb C})$ be a smooth map. Then 
$g$ can be viewed as an automorphism of the trivial
complex vector bundle ${\Bbb C}^N\rightarrow X$ over $X$.
Thus
for any $b\in B$, $g$ induces a bounded map $g_b$ from
$L^2(S(TZ_b)\otimes E_b\otimes {\Bbb C}^N)$ to itself by acting as an identity on
$L^2(S(TZ_b)\otimes E_b)$.

Also, for the spectral section $P$, $P_b$ induces a map on 
$P_bL^2(S(TZ_b)\otimes E_b\otimes {\Bbb C}^N)$ by acting as identity on ${\Bbb C}^N$, denoted again by $P_b$.
 Then the Toeplitz family $T_g=\{ T_{g,b}\}_{b\in B}$ is a family of
Toeplitz operators defined by
$$T_{g,b}=P_bg_bP_b:P_bL^2(S(TZ_b)\otimes E_b\otimes {\Bbb C}^N)\rightarrow
P_bL^2(S(TZ_b)\otimes E_b\otimes {\Bbb C}^N).$$

The Chern character of its index bundle is expressed in terms of topological data as follows.

\begin{theo} \label{toe}The following identity holds,
$${\rm ch}({\rm ind}\, T_g)=-\int_{Z}\hat{A}(TZ){\rm ch}(E){\rm ch}(g)\ \ 
{\rm in}\ \ H^*(B),  $$
where ${\rm ch}(g)$ is the odd Chern character associated to $g$ (cf. [AH]).
\end{theo}

Just as the index of Toeplitz operators can be expressed via spectral flow, we interpret the 
index bundle of a Toeplitz family via higher spectral flow. This generalizes a result of 
Booss-Wojciechowski \cite[Theorem 17.17]{bw}.

\begin{theo} We have the following identity in $K(B)$,
$$ {\rm ind}\, T_g ={\rm sf}\{D^E, gD^Eg^{-1}\}.$$
\end{theo}

  A particularly nice application of this result is a local version of Theorem~\ref{toe}, which 
can be viewed as a odd analogue of the Bismut local index theorem for a family of Dirac 
operators on even dimensional manifolds. 

We will also use higher spectral flow to prove a generalization of the family index theorem 
for manifolds with boundary  \cite{bc2}, \cite{bc3}, \cite{mp}.

Let us mention that Douglas-Kaminker [DK] has also proposed a kind of higher
spectral flow concept independently. On the other hand, Wu [W] has extended our
concept of higher spectral flow to noncommutative geometry.

This paper is organized as follows. In Section 1, we introduce the concept of
higher spectral flow, and show that it is a generalization of the classical
concept of spectral flow inroduced in [APS1]. In Section 2 we explore the
relationship between the Chern character of the higher spectral flow and the
$\hat{\eta}$-forms. In Section 3, we consider  periodic families and prove
Theorem 0.2. In Section 4, we introduce the concept of Toeplitz family, give a proof 
of Theorem 0.5, as well as a heat kernel proof of Theorem 0.4. Finally in 
Section 5 we prove two results for family indices of manifolds with boundary, which
generalize the results in [DZ1], [MP] respectively.

The results of this paper have been announced in [DZ3].

$\ $

{\bf Acknowledgements.} We would like to thank K. Wojciechowski for a helpful
discussion concerning Toeplitz operators.
Part of this work is done while the first author was 
visiting the Nankai Institute of Mathematics. He is very grateful to the Chinese National 
Science Foundation for partial financial support and the Nankai Institute for its hospitality.  
The second author would like to thank the Courant Institute for its hospitality. 

\sect{ Spectral sections and higher spectral flows}

In this section we will first recall the notion of spectral flow and its
basic properties in a). Then we will give, in b),  a new interpretation using the notion of spectral 
sections. This new interpretation is the basis for the higher spectral
flow. At the end of the section, in c),  we extend this interpretation using a generalized notion of spectral section. This gives us the flexibility we need in proving Theorem 4.4.

$\ $

{\bf a). Spectral flow}

$\ $

 We take the definition of spectral flow as in \cite{aps}. Thus, if $D_s, \ 0\leq s \leq 1$, is 
a curve of self-adjoint Fredholm operators, the spectral flow ${\rm sf}\{D_s\}$ counts the net 
number of eigenvalues of $D_s$ which change sign when $s$ varies from $0$ to $1$. 
(Throughout the paper a family always means a {\em continuous} (and sometimes
{\it smooth}) family, and a curve 
always means a one parameter family.) The following proposition collects some of its most 
basic properties from [APS1] and [APS2].

$\ $

{\bf Proposition 1.1}. {\it 1). If $D_s, \ 0\leq s \leq 1$, is a curve of self-adjoint Fredholm operators, 
and $\tau \in [0, 1]$, then
$$ {\rm sf}\{D_s, [0, 1]\} = {\rm sf}\{D_s, [0, \tau]\} + {\rm sf}\{D_s, [\tau, 1]\}. 
\eqno (1.1)$$
2). If $D_s, \ 0\leq s \leq 1$, is a smooth curve of self-adjoint elliptic pseudodifferential operators 
on a closed manifold, and $\bar{\eta}(D_s) =\hf(\eta(D_s)+\dim \ker D_s)$
is the reduced $\eta$ invariant of $D_s$ in the sense of
Atiyah-Patodi-Singer [APS2], then $\bar{\eta}$ is ${\rm mod}\ {\Bbb Z}$ smooth and
$$ {\rm sf}\{D_s\}= -\int^1_0 \frac{d\bar{\eta}(D_s)}{ds} \, ds + \bar{\eta}(D_1) -\bar{\eta}(D_0) . \eqno (1.2)$$
3).  If $D_s, \ 0\leq s \leq 1$, is a periodic one parameter family of self-adjoint Dirac type 
operators on a closed manifold, and $\tilde{D}$ is the corresponding Dirac type operator 
on the mapping torus, then
$$ {\rm sf}\{D_s\}= \ind\ \tilde{D} . \eqno (1.3)$$}

$\ $

{\bf b). Spectral sections}

$\ $

The notion of spectral section is quite new and is introduced by Melrose-Piazza [MP] in their study of the family 
index of Dirac operators on manifolds with boundary. It can be viewed as a generalized 
Atiyah-Patodi-Singer boundary condition and can be defined for a family of  self adjoint 
first order elliptic pseudodifferential  operators.
\newline

{\bf Definition 1.2}. A spectral section for a family of
first order elliptic pseudodifferential  operators $D=\{ D_z \}_{z\in B}$ is a family of self adjoint 
pseudodifferential projections $P_z$  on the $L^2$-completion of the domain of
$D_z$ such that for some smooth function $R:\, B 
\longrightarrow {\Bbb R}$ and every $z\in B$
$$ D_zu=\lambda u \Rightarrow \left\{ \begin{array}{l} P_zu=u \ \mbox{if}\ \lambda > 
R(z) \\ P_zu=0 \ \mbox{if}\ \lambda < - R(z).  \end{array} \right.$$ 

We list in the following the basic properties of spectral section. The reader is referred to 
\cite{mp} for detail.

$\ $

{\bf Proposition 1.3}. {\it Let $D=\{ D_z \}_{z\in B}$ be a family of
first order elliptic pseudodifferential  operators, and assume that the parameter space $B$ is compact. Then \\
A). There exists  a spectral section for $D$ if and only if the (analytic) index [AS2] of the family vanishes: $\ind \, D =0 $ in 
$K^1(B)$. In particular,  spectral section for a one parameter family always exists. \\
B). Given spectral sections $P$, $Q$, there exists a spectral section $R$ such that
$PR=R$ and $QR=R$. Such a spectral section will be called a majorizing spectral section. \\
C). If $R$ majorizes $P$: $PR=R$, then $\ker \{ P_bR_b:{\rm Im}(R_b)\rightarrow
{\rm Im}(P_b)\}_{b\in B}$  forms a vector bundle on $B$, denoted by
$[R-P]$. Hence for any two spectral sections $P$, $Q$, the difference element $[P-Q]$ 
can defined as an element in $K(B)$, as follows:
$$ [P-Q] = [R-Q] - [R-P] , \eqno (1.5)$$
for any majorizing spectral section $R$.  \\
D). If $P_1$, $P_2$, $P_3$ are spectral sections, then
$$ [P_3-P_1] =[P_3-P_2] + [P_2-P_1] . \eqno (1.6)$$
E). The $K$-group $K(B)$ is generated by all these difference elements. }

$\ $

Now let $D_s$ be a curve of self adjoint elliptic pseudodifferential operators. Let $Q_s$ be 
the spectral projection onto the direct sum of eigenspaces of $D_s$ with nonnegative 
eigenvalues (the APS projection). The following theorem provides a link between the 
above two notions.

$\ $

{\bf Theorem 1.4}. {\it Let $P_s$ be a spectral section of $D_s$. Then $[Q_1-P_1]$ 
defines an 
element of $K^0(pt) \cong {\Bbb Z}$ and so does $[Q_0-P_0]$. Moreover 
the difference $[Q_1-P_1] -  [Q_0-P_0]$ is independent of the choice of the spectral 
section $P_s$, and it computes the spectral flow of $D_s$:
$$ {\rm sf}
\{D_s\}=[Q_1-P_1] -  [Q_0-P_0]. 
\eqno (1.7)$$}

{\it Proof}.  We first show the independence of the choice of spectral section. Thus let $P$, $Q^\prime$ be 
two spectral sections and $R$ be a majorizing spectral section. Since $[R-P]$ defines a 
vector bundle on $[0, 1]$, we have in particualr $\dim [R_1 -P_1]=\dim [R_0-P_0]$,  and 
similarly, $\dim [R_1 -Q_1^\prime]=\dim [R_0-Q_0^\prime]$. This, together with the basic property D), 
gives us the independence. For (1.7), we note that both sides are additive with respect 
to the subdivision of the interval. Thus, by dividing into sufficiently small intervals if necessary, we can 
assume that there exists  a positive real number $a$ with $\pm a \not\in \spec(D_s)$ for all $s$. Since we are free to choose any spectral section to compute the right hand side of (1.7), we choose 
$P_s$ to be the orthogonal projection onto the direct sum of eigenspaces of $D_s$ 
with eigenvalues greater than $a$. This reduces the problem to the finite dimensional case, 
where it can be easily verified. 
\qed

$\ $

Theorem 1.4 leads us to the notion of higher spectral flow, when we consider higher dimensional 
families.

Let $\pi:\,  X \longrightarrow B$ be a smooth fibration with the typical fiber $Z$ an odd 
dimensional closed manifold and $B$ compact. As we mentioned before, a family of self adjoint elliptic 
pseudodifferential operators on $Z$, parametrized by $B$, will be called a $B$-family. 
Consider
a curve of $B$-families, $D_u=\{D_{b,u}\}$, $u\in [0,1]$. 

Assuming that the index bundle of $D_0$ vanishes, the homotopy invariance of the index 
bundle
then implies that the index bundle of each $D_u$ vanishes.  Let $Q_0$, $Q_1$ be spectral 
sections of $D_0$, $D_1$ respectively. If we consider the total family $\tilde{D}=\{D_{b, 
u}\}$ parametrized by $B\times I$, then there is a total spectral section 
$\tilde{P}=\{P_{b,u}\}$. Let $P_u$ be the restriction of $\tilde{P}$ over $B\times 
\{u\}$.
\newline

{\bf Definition 1.5}.  The (higher) spectral flow ${\rm sf}\{(D_0, Q_0),\ (D_1, Q_1)\}$ between the pairs $(D_0, 
Q_0)$, $(D_1, Q_1)$ is an element in $K(B)$ defined by
$$ {\rm sf}\{(D_0, Q_0),\ (D_1, Q_1)\} = [Q_1-P_1] -  [Q_0-P_0] \in K(B). \eqno (1.8)$$

The definition is independent of the choice of the (total) spectral section $\tilde{P}$, as it follows 
again from the basic properties C), D).
When $D_u$, $u\in S^1$ is a periodic family, we choose $Q_1=Q_0$. In this case the 
spectral flow turns out to be independent of $Q_0=Q_1$ and therefore defines an 
intrinsic invariant 
of the family, which we denote by ${\rm sf}\{D_u\}$.
\newline

Although this notion is defined generally, for the most part of this paper we are going to restrict our attention to
$B$-families of Dirac type operators, defined as follows.

 For simplicity we make the 
assumption that the vertical bundle $TZ\longrightarrow X$ is spin and we fix a spin 
structure\footnote{Our discussion extends without difficulty to the 
more general case when there are smoothly varying ${\Bbb Z}_2$-graded Hermitian 
Clifford modules over the fibers, with graded unitary connections.}.
Let $g^{TZ}$ be a metric on $TZ$. We use $S(TZ)$ to denote the spinor bundle of $TZ$. 
Let $(E, g^E)$ be a Hermitian vector bundle over $X$, and $\nabla ^E$ a Hermitian
connection on $E$.

For any $b\in B$, one has a canonically defined self-adjoint (twisted) Dirac operator
$$ D_{ b}^E :\ \Gamma((S(TZ)\otimes E)|_{Z_b}) \longrightarrow  
\Gamma((S(TZ)\otimes E)|_{Z_b}). 
\eqno (1.9)$$
This defines a smooth family of (standard, twisted) Dirac operators $D^E$ over $B$.
\newline

{\bf Definition 1.6}.  By a $B$-family of Dirac {\em type} operators on $(X, E)$, we mean a smooth 
family of first order self adjoint differential operators $D=\{D_b\}_{b\in B}$ whose principal 
symbol is given by that of $\{D_{ b}^E\}_{b \in B}$.

$\ $

Thus, for any $B$-family of Dirac type operators $\{D_b\}_{b\in B}$, there is a natural 
homotopy between $\{D_b\}_{b\in B}$ and $\{D_{ b}^E\}_{b \in B}$ through $B$-
families of Dirac type operators given by
$$
\tilde{D}_b(t)= (1-t) D_b + t D_{ b}^E. 
\eqno (1.10)$$

Since the definition of higher spectral flow requires the existence of a spectral section, we 
now make the following
\newline

{\bf Basic Assumption 1.7.} We assume that the (canonical) family $D^E$ has vanishing 
index bundle.
\newline 

A typical example satisfying our basic assumption is the family of signature operators. 
More generally a $B$-family whose kernels have constant dimension
will always satisfy the basic assumption. Another class of examples comes from the 
boundary family of a family of Dirac type operators on even dimensional manifolds with boundary.

 Now we can speak of the higher spectral flow of a curve of $B$-families of Dirac type 
operators, given the basic assumption.
\newline

{\bf c). Generalized spectral sections}

$\ $

The rest of this section is devoted to an extension of the notion of spectral section. This is what 
we call generalized spectral section. The higher spectral flow can also be defined in terms of generalized spectral sections. This will come into play later when we discuss Toeplitz families. 
\newline

{\bf Definition 1.8} ([DZ2]).  A generalized spectral section $Q$ of  a family of  self adjoint 
first order elliptic pseudodifferential  operators $D=\{ D_z\}_{z\in B}$  is a continuous
family of self-adjoint zeroth order pseudodifferential projections $Q=\{ Q_z\}
_{z\in B}$ whose principal symbol is the same
as that of a (in fact any) spectral section $P$.

$\ $

{\bf Remark 1.9}.  There are actually  many equivalent definitions of generalized
spectral sections. For example one can use the language of
infinite dimensional Grassmannian
discussed in \cite{bw}. Also one can use the infinite dimensional Lagrangian subspaces
as used by Nicolaescu in [N]. 
One of the most important examples of  generalized spectral sections is the Cald\'eron 
projection (cf. [BW]). It has been used in an essential way in our proof of the
splitting formula for family indices [DZ2].

$\ $

As we have seen, the crucial ingredient in the definition of higher spectral flow is the 
difference element of two spetral sections. This notion generalizes
to generalized spectral sections.

Thus let $Q_1$, $Q_2$ be two generalized spectral sections of $D$.
For any $z\in B$, set
$$T_z(Q_1,Q_2)=Q_{2,z}Q_{1,z}:{\rm Im}(Q_1)\rightarrow
{\rm Im}(Q_2).\eqno (1.11)$$
Then $T(Q_1,Q_2)=\{ T_z(Q_1,Q_2)\}_{z\in B}$ defines a continuous family of 
Fredholm
operators over $B$. Thus according to Atiyah-Singer [AS1], it determines an element
$$[Q_1-Q_2]=\mbox{ind}\, T(Q_1,Q_2)\in K(B).\eqno (1.12)$$
In the special case where $Q_1$ and $Q_2$ are two 
spectral sections, 
one verifies easily that $[Q_1-Q_2]$ is the same as the difference element defined
by Melrose-Piazza \cite{mp}.
\newline

Two generalized spectral sections $Q_1$, $Q_2$ of $D$ are
said to be homotopic to each other if there is a continuous curve of generalized
spectral sections $P_u$, $0\leq u\leq 1$, of $D$ such that $Q_1=P_0$, $Q_2=P_1$.

$\ $

{\bf Proposition 1.10}. {\it  1).  If $Q_i$, $i=1,2,3$, are three generalized spectral sections 
of $D$ such that $Q_1$ and $Q_2$ are homotopic to each other, then one has
$$[Q_1-Q_3]=[Q_2-Q_3] \ \ \mbox{in} \ \ K(B). \eqno (1.13)$$
2).  If $Q_i$, $i=1,2,3$, are three generalized spectral sections
of $D$, then the following identity holds in $K(B)$,
$$[Q_1-Q_2]+[Q_2-Q_3]=[Q_1-Q_3].\eqno (1.14)$$}

\Pf. These follow easily from some elementary arguments concerning
Fredholm families. \qed

The generalized spectral sections are more flexible than the spectral sections since they only have to have the right symbol. On the other hand, we still have

$\ $

{\bf Theorem 1.11.} {\it The higher spectral flow (1.8) can be computed using generalized spectral sections. That is, we can choose $\tP$ to be a generalized spectral section of $\tD$.}

$\ $

{\it Proof.} The point here is that the right hand side of (1.8) is also independent of the choice of generalized spectral sections. To see this, let $\tP$ and $\tP'$ be two generalized spectral sections of the total family $\tD$ parametrized by $B\times I$. Then $T(P_u, P_u')$ defines a curve of Fredholm families. By the homotopy invariance of the family index, we have
\[ [P_1-P_1']=\ind T(P_1, P_1') = \ind T(P_0, P_0')= [P_0-P_0'].\]
\qed

$\ $

{\bf Remark 1.12.}  Given generalized spectral sections $P$, $Q$, there exists  a generalized spectral section $R$ such that
$$(1-Q)L^2(Z,S(TZ)\otimes E|_Z)\subset (1-R)L^2(Z,S(TZ)\otimes E|_Z)$$
and that the family of Fredholm operators $T(1-R,1-P)$ has vanishing cokernels,
which implies that the family $T(P,R)$ has vanishing cokernels.
This can be obtained by applying the procedure in [AS1] to the family of
Fredholm operators
$$T(1-Q,1-P):(1-Q)L^2(Z,S(TZ)\otimes E|_Z)\rightarrow (1-P)L^2(Z,S(TZ)\otimes E|_Z). $$
The generalized spectral section $R$ is an analogue of a majorizing spectral section.

\sect{ $\heta$-form and the Chern character of higher spectral flows}

In this section we show that the Chern character of  the higher spectral flow can be 
calculated by the heat kernel methods through $\heta$-forms. For simplicity 
we will restrict ourselves to Dirac type operators while at the same time we 
allow the freedom to use the superconnections generalizing the Bismut superconnection 
\cite{b}. In this sense the $\heta$-form of this section can be seen as certain generalizations 
of those of Bismut-Cheeger \cite{bc1} and Melrose-Piazza \cite{mp}.

This section is organized as follows. In a), we 
recall some basic results of Melrose-Piazza~\cite{mp}. In b), we introduce a generalized 
superconnection. In c), we define the $\heta$-forms and prove certain basic properties. In 
d), we prove a relative formula generalizing a result of ~\cite{mp}. Finally in e), we prove 
the main result of this section which expresses the Chern character of higher spectral 
flow via $\heta$-forms.
\newline

{\bf a). Melrose-Piazza operator of a spectral section }
\newline

Now, let $D=\{D_b\}_{b\in B}$ be a $B$-family of self adjoint Dirac type operators as defined 
in Section 1, and let $P=\{P_b\}_{b\in B}$ be a spectral section of $\{D_b\}_{b\in B}$. 
Then by~\cite{mp} there exists a family of zeroth order finite rank pseudodifferential 
operators $\{A_b\}_{b\in B}$, 
$$ A_b :\ \Gamma((S(TZ)\otimes E)|_{Z_b}) \longrightarrow  \Gamma((S(TZ)\otimes 
E)|_{Z_b}),$$ 
such that, for any $b\in B$, \\
(i). $\tilde{D}_b = D_b + A_b$ is invertible; \\
(ii). $P_b$ is precisely the Atiyah-Patodi-Singer projection [APS2] of $\tilde{D}_b$.
\newline

{\bf Definition 2.1}.  We call this $A=\{A_b\}_{b\in B}$ a Melrose-Piazza operator associated to the 
spectral section $P$.
\newline

{\bf b). Superconnections associated to Dirac type families}
\newline

We now choose a connection for the fibration which  amounts to a splitting 
\be 
TX = TZ \oplus T^HX.
\ee
We also have the identification $T^HX=\pi^*TB$.  

Endow $B$ with a metric $g^{TB}$ and let $g^{TX}$ be the metric defined by
\be
g^{TX} = g^{TZ} \oplus \pi^*g^{TB}.
\ee
Let $P$, $P^{\perp}$ be the orthogonal projections of $TX$ onto $TZ$, $T^HX$ 
respectively and denote by $\nabla^{TX}$, $\nabla^{TB}$ the Levi-Civita connections of 
$g^{TX},\, g^{TB}$ respectively. 
Following Bismut \cite{b}, let $\nabla^{TZ}$ be the connection on the vertical bundle 
defined by
$\nabla^{TZ}=P\nabla^{TX}P$. This is a connection compatible with the metric 
$g^{TZ}$ and is independent of the choice of the metric $g^{TB}$.

Then the connection lifts to a connection on the spinor bundle. Also following Bismut we 
view $\Gamma(S(TZ) \otimes E)$ as the space of sections of an infinite dimensional vector 
bundle $H_{\infty}$ over $B$, with fiber 
\be
H_{\infty, b} = \Gamma(S(TZ_b)\otimes E_b).
\ee
Then $\nabla^{S(TZ)\otimes E}$ determines a connection on $H_{\infty}$  by
the prescription:
\be
\tilde{\nabla}_X h = \nabla_{X^H} h,
\ee
where $Y^H\in \Gamma (T^HX)$ is the horizontal lift of $Y\in \Gamma (TB)$.

Now for any $1\leq i \leq \dim B$, let  
$$B_i \in \Omega^i(T^*B)\hat{\otimes} 
\Gamma(\End(S(TZ){\otimes} E)) = \Omega^i(T^*B)\hat{\otimes}{\rm cl}(TZ) \otimes \End(E)$$ 
be an 
odd element (with respect to the natural ${\Bbb Z}_2$-grading). And let $A$ be a 
Melrose-Piazza operator associated to the spectral section $P$. Furthermore, choose
a cut-off function $\rho$ in $t$: $\rho(t)=1$ when $t>8$ and $\rho(t)=0$ when 
$t<2$.
\newline

{\bf Definition 2.2}.  For any $t>0$, the (rescaled) superconnection $B_t$ is a superconnection on 
$H_{\infty}$ given by
\be
B_t = \tilde{\nabla} + \sqrt{t}(D+ \rho(t) A) + \sum^{\dim B}_{i=1} t^{\o{1-i}{2}}B_i, 
\ee
We also set $B=B_1= \tilde{\nabla} + D + \sum^{\dim B}_{i=1} B_i$. (This should not be confused with the base manifold.)
\newline 

{\bf c). $\heta$-functions and their residues}
\newline

We define, for $\Re(s) \gg 0$, the $\heta$-function for the superconnection $B_t$ as
follows
\be \heta(B, A, s)=\o{1}{\sqrt{\pi}} \int^{\infty}_0 t^{\o{s}{2}} 
\Tr^{\rm even}[\o{dB_t}{dt} e^{-B_t^2}] dt.  \ee

This defines an even differential form on $B$ which depends holomorphically on $s$ for 
$\Re(s) \gg 0$.
By  standard arguments it extends to a meromorphic function of $s$ in the whole 
complex plane with only simple poles. 

$\ $

{\bf Remark 2.3}.  Note that our definition has an extra factor of $\hf$, coming from
${dB_t\over dt}$,
in comparing with those of  \cite{bc1}, \cite{mp}. 

$\ $

{\bf Proposition 2.4}. {\it The residue of $\heta(B, A, s)$ at $s=0$ is an exact form.}

{\it Proof}.  Our proof consists of two steps. We first show that, modulo exact forms, 
$\mbox{Res}_{s=0}\{\heta(B, A, s)\}$ is invariant under the smooth deformation of Dirac 
type families. The main ingredient here is the transgrssion formula (2.8).

But note first of all, since $\rho(t)=0$ for $t\leq 1$, $\mbox{Res}_{s=0}\{\heta(B, A, s)\}$
does not depend on either $A$ or $P$. 

Now let $\{D_u\}_{u\in[0, 1]} $ be a smooth curve of $B$-families $D_u=\{D_{b, u}\}_{b\in 
B}$. We view $\tilde{D}=\{D_{b, u}\}_{b\in B, u\in [0, 1]}$ as a $B\times I$-family 
of Dirac type operators. Note that this family also verifies the basic assumption.

Let $\tilde{P}$ be a spectral section of $\tilde{D}$ and $\tilde{A}$ a Melrose-Piazza operator for 
$\tilde{P}$. By pulling back $B_i$ to an odd element  $\tilde{B}_i$ on $B\times I$, we 
obtain a superconnection $\tilde{B}_t$ on $B\times I$:
\be
\tilde{B}_t =\o{\pt}{\pt u} du + \tilde{\nabla} + \sqrt{t}(\tilde{D}+ \rho(t) \tilde{A}) + 
\sum^{\dim B}_{i=1} t^{\o{1-i}{2}}\tilde{B}_i. 
\ee

Proceeding as in \cite{q}, \cite{b}, we deduce
\be
\o{\pt}{\pt t} \Tr^{\rm odd} [ \exp (-\tilde{B}_t^2)] = (\o{\pt}{\pt u} du + 
d_B)\Tr^{\rm even}[\o{\pt \tBt}{\pt t} \exp(-\tBt^2)].
\ee
Therefore, for  $\Re(s)\gg 0$, 
\be   
(\o{\pt}{\pt u} du + d_B)\int^{\infty}_0 t^{\o{s}{2}}\Tr^{\rm even}[\o{\pt \tBt}{\pt t} \exp(-
\tBt^2)] dt =-\o{s}{2} \int^{\infty}_0 t^{\o{s}{2}-1} \Tr^{\rm odd} [ \exp (-\tilde{B}_t^2)] 
dt,
\ee
where we have used  the fact that $\tD + \tA$ is invertible.

Now when $\dim Z$ is odd, by proceeding as in \cite{b}, \cite{bf}, one has the 
following small time asymptotics
\be
\Tr^{\rm odd} [ \exp (-\tBt^2)] = \frac{a_{-k}}{t^k} +{ a_{-k+1}\over {t^{k-1}}} + \cdots + a_0 + 
O(t)
\ee
for some $k \in {\Bbb Z}$ as $t \ra 0$.
This gives
\be  
\Res_{s=0} \left\{ -\frac s2 \int_0^\infty t^{\frac s2 -1} \Tr^{\rm odd} [\exp (-\tBt^2)] dt 
\right\} =0.
\ee
For any $u \in [0,1]$, set
\be
B(t,u) = \tilde{\nabla} + \sqrt{t} (D_{b,u} + \rho (t) A_{t,u}) + \sum_{i=1}^{\dim B} 
t^{\frac{1-i}{2}} B_{i,u}.
\ee
Then
$$\tBt ^2  =  (\frac{\partial}{\partial u} du + B(t,u))^2 \\
 =  B(t,u)^2 + du \frac{\partial B(t,u)}{\partial u}.\eqno (2.13)$$
Thus
$$\Tr^{\rm even} \left[ \frac{\partial \tBt}{\partial t} \exp (-\tBt^2)\right]  =  
\Tr^{\rm even} \left[ \frac{\partial B(t,u)}{\partial t} \exp (-B(t,u)^2 - du \frac{\partial 
B(t,u)}{\partial u} ) \right] $$ 
$$ =   \Tr^{\rm even} \left[ \frac{\partial B(t,u)}{\partial t} \exp (-B(t,u)^2) \right] $$ 
$$+ \frac{\partial}{\partial \tau} \left\{ \Tr^{\rm even} \left[ \frac{\partial B(t, u)}{\partial t} 
\exp (-B(t,u)^2 - \tau \frac{\partial B(t,u)}{\partial u} du) \right] \right\} |_{\tau=0}.\eqno (2.14)$$
 
 From  (2.14), (2.11) (2.9) and (2.6), one deduces that
$$
\Res_{s=0} \left\{ \frac{\partial \heta (B_u,A_u,s)}{\partial u} \right\} du + $$
$$  d_B\left\{ 
\frac{1}{\sqrt{\pi}} \Res_{s=0} \{ \int_0^\infty t^{\frac s2} \frac{\pt}{\pt \tau} \{ 
\Tr^{\rm even} [ \frac{\partial B(t,u)}{\partial t} \exp (-B(t,u)^2 - \tau  \frac{\partial B(t,u)}{\partial u} du)]\}\left|_{\tau =0} \right. \} \right\}$$
$$ =0.\eqno (2.15)$$

Therefore
$$\Res_{s=0} \left\{ \frac{\partial \heta (B_u,A_u,s)}{\partial u} \right\} \in d\Omega (B).\eqno (2.16)$$

With this invariance property at  hand, we can now easily finish our proof. Let  $D^E$
be the $B$-family of (standard) Dirac operators, and $B^E$ the Bismut superconnection.
Then by~\cite{bc1}, \cite{mp}, 
$$ \mbox{Res}_{s=0}  \{ \heta(B^E, A_1, s)\}=0.\eqno (2.17)$$
Applying the invariance property to the curve of $B$-families
$ D_u=(1-u)D+uD^E$ yields the desired result. \qed
\newline

{\bf d).  $\heta$-forms}
\newline

With the help of  Proposition 2.4  we can now define the $\heta$-forms.

$\ $

{\bf Definition 2.5}.  The $\heta$-form (associated to $B$ and $A$) is defined as
$$ \heta(B, A)=\left\{ \heta(B, A, s)-\o{\mbox{Res}_{s=0}\{\heta(B, A, s)\}}{s} 
\right\}_{s=0}. \eqno (2.18)$$ 

$\ $

{\bf Remark 2.6}.  The $\heta$-form for a general superconnection is defined in \cite{bc1} in the  finite 
dimensional case.
\newline

A variational argument shows

$\ $

{\bf Proposition 2.7}. {\it The value of $\heta(B, A)$ in $\Omega^*(B)/d\Omega^*(B)$ is independent 
of the choice of the cut-off function $\rho$ and the Melrose-Piazza operator $A$.}

{\it Proof}. If $\gamma (t)$ is another cut-off function, we will set 
$\rho_u(t)=(1-u)\rho(t) +u\gamma (t). $

Also, if $A^\prime$ is another Melrose-Piazza operators, by \cite{mp}, there is a smooth family 
$A(u)$ of Melrose-Piazza operators such that $A(0)=A$, $A(1)=A^\prime$.

Then the rest of the proof follows along the same line as  that in the proof of Proposition
2.4. Namely, the variation of the $\heta$-forms is, up to an exact form, given by a local 
term (Cf. formula (2.24) below). This local term is the constant term in the asymptotic expansion (2.22). Since $\rho(t)=0$ for $t\leq 2 $ (and all $u$), and $D$ is independent of $u$,  this term is zero. \qed

$\ $

This enables us to give an intrinsic definition.
\newline

{\bf Definition 2.8}. Given a spectral section $P$ and a superconnection $B$, the $\heta$-form $\heta(D, 
P)$ is defined to be an element of $\Omega^*(B)/d\Omega^*(B)$ determined by $\heta(B, 
A)$ for any choice of $A$. 
\newline

{\bf e). Difference element and $\heta$-forms}
\newline

The dependence of $\heta$-form on the spectral section is also well understood.
The following is a slight generalization of a result of \cite{mp}.

$\ $

{\bf Thereom 2.9}. {\it If $P_0$ and $P_1$ are spectral sections of the family $D$, then
$${\rm ch}(P_1-P_0) = \heta(D, P_1)-\heta(D, P_0) \ {\rm in} \ H^*(B).
\eqno (2.19)$$}

{\it Proof}. We note first that if one fix the family $D$, the spectral section $P$, and varies the other terms in the superconnection, the variation of the $\heta$-form will be given by a local term, which is independent of the spectral section. This shows that the right hand side of (2.19) is independent of such variations. Now we can connect any superconnection to the Bismut superconnection through such a deformation. Since (2.19) is true for the Bismut superconnection by [MP], the proof is complete. 
\qed
\newline

{\bf f). $\heta$-forms and the Chern character of higher spectral flow}
\newline

Now consider a smooth curve of $B$-families of Dirac type operators $\tD = \{ 
\tD_u\}_{u\in [0,1]}$. We view it as a total family over $B\times I$. Also, let $\tP$ be a 
total spectral section for $\tD$ and $\tA$ a corresponding Melrose-Piazza operator. As 
before, set
$$ B(t,u) = \tilde{\nabla} + \sqrt{t} (D_{u} + \rho (t) A_{u}) + \sum_{i=1}^{\dim B} 
t^{\frac{1-i}{2}} B_{i,u}. \eqno (2.21)$$

Proceeding as in \cite{b}, \cite{bf}, one verifies the following asymptotic expansion
$$\o{\pt}{\pt \tau} \left\{ \Tr^{\rm even} [ \exp (-B(t,u)^2-\tau \o{\pt B(t, u)}{\pt 
u})]\right\}_{\tau=0} = 
\frac{a_{-k}}{t^k} + {a_{-k+1}\over t^{k-1}} + \cdots + a_0 + O(t)\eqno (2.22)$$
for some $k \in {\Bbb Z}$ as $t \ra 0$. 
Here $a_i=a_i(u)$ are smooth forms on $B$.

$\ $

{\bf Theorem 2.10}. {\it Let $P_0$, $P_1$ be the restrictions of $\tP$ at $u=0, \, 1$ 
respectively. Then
$$\heta(B_1, P_1) - \heta(B_0, P_0) = \o{1}{\sqrt{\pi}}\int^1_0 a_0(u) du, \ {\rm in} \ 
\Omega^*(B)/d\Omega^*(B). \eqno (2.23)$$ }

{\it Proof}.  From (2.13), (2.9), (2.14), one deduces
$$\o{\pt \heta(B_u, A_u, s)}{\pt u} du + d_B \o{1}{\sqrt{\pi}} \int^{\infty}_0 
t^{\o{s}{2}} \o{\pt}{\pt \tau} \{ \Tr^{\rm even}[\o{\pt B_u(t)}{\pt t} \exp(-B_u(t)^2+\tau 
\o{\pt B_u(t)}{\pt t} du)] \}_{\tau=0} dt$$
$$ = -\o{s}{2\sqrt{\pi}}  \o{\pt}{\pt \tau} \{ \int^{\infty}_0 t^{\o{s}{2}-1} 
\Tr^{\rm odd}[\exp(-B_u(t)^2-\tau \o{\pt B_u(t)}{\pt u} du)] dt \}_{\tau=0}.
\eqno (2.24)$$

We now let $s\ra 0$ and equate the constant terms in the corresponding asymptotic 
expansions by standard methods, and then integrate from 0 to 1.  \qed

$\ $

{\bf Remark 2.11}. The main point here is that the $\heta$-forms $\heta(B_u, P_u)$ is continuous in 
$u$. Thus, there is no ``jump phenomena".

$\ $

The following theorem generalizes the well known relationship between the spectral flow 
and the $\eta$ invariant ([APS1]).

$\ $

{\bf Theorem 2.12}. {\it Let $D_u$ be a curve of $B$-families of Dirac type operators and $Q_0$, $Q_1$ spectral 
sections of $D_0$, $D_1$ respectively. Let $\tB_t(u) = \tilde{\nabla} +\sqrt{t} (D_u + 
\rho(t) \tilde{A}) +\sum_{i\geq 1} t^{\o{1-i}{2}} B_i(u)$ be a curve of superconnections. 
Then
we have the following identity in $H^*(B)$:
$${\rm ch}({\rm sf}\{(D_0, Q_0),\ (D_1, Q_1)\}) = \heta(D_1, Q_1) -\heta(D_0, Q_0)
  -\o{1}{\sqrt{\pi}}  \int^1_0 a_0(u) du ,\eqno (2.25)$$ 
where $a_0(u)$ is determined by (2.22).}

{\it Proof}.  This follows  from the definition of higher spectral flow, Theorems 
2.9 and 2.10. \qed

$\ $

Theorem 2.12 provides a way to compute, analytically, the Chern character of higher 
spectral flow. We will use formula (2.25) to compute the higher spectral flow of both 
a periodic family and a Toeplitz family.

\sect{ Periodic family index and the higher spectral flow}

In this section, we use the results of Section 2 to calculate the Chern character
of the higher spectral flow of periodic families
and show that it equals the Chern character of the family index 
of the associated families of Dirac operators on $even$ dimensional manifolds.

This section is organized as follows. In a), we introduce the periodic family
and the associated higher spectral flow. In b), we examine the 
relation between the Chern character of the higher spectral flow and of the 
family index.

$\ $

{\bf a). Periodic families}

$\ $

By definition, a periodic  $B$-family of self-adjoint Dirac type 
operators is just a closed curve of $B$-families of operators. Equivalently, one
can also view it as a $u\in [0,1]$ family of $B$-families of Dirac type
operators, such that $D_0=D_1$.

In this special case, if we give $D_0,\ D_1$ the same spectral section $Q$ and let
$\tilde{P}$ be a total spectral section over $B\times I$, then one has
$${\rm sf} \{ (D_0,Q),(D_1,Q)\} =[P_1-Q]-[P_0-Q]$$
$$=[P_1-P_0].\eqno (3.1)$$

Thus one has

$\ $

{\bf Proposition 3.1}. {\it The higher spectral flow ${\rm sf}\{ (D_0,Q),(D_1,Q)\}$
for  the periodic family does not depend on $Q$.}

$\ $

Therefore ${\rm sf}\{ (D_0,Q),(D_1,Q)\}$ defines an intrinsic invariant of the
periodic family. We will denote it by ${\rm sf}\{ D\}$.

$\ $

{\bf b). Relations of higher spectral flows with index bundles for periodic
families}

$\ $

Now let $D=\{ D_u\} _{u\in S^1}$ be a periodic family of $B$-families of
(self-adjoint) Dirac type operators $\{ D_{b,u}\} _{b\in B,u\in S^1}$.
Then for any $b\in B$, $\{ D_{b,u}\}_{u\in S^1}$ defines a natural Dirac
operator $D_b^\prime$ on $S^1\times Z_b$:
$$D_b^\prime = \left(
 \begin{array}{cc} 0 &  \frac{\pt}{\pt u} + D_{b,u} \\
- \frac{\pt}{\pt u} + D_{b,u} & 0
\end{array} \right) : 
$$
$$\Gamma ({\Bbb C}^2\otimes S(TZ_b)\otimes E_b)\rightarrow \Gamma
({\Bbb C}^2\otimes S(TZ_b)\otimes E_b),\eqno (3.2)$$
where we have lift $TZ_b,\ E_b$ canonically to $S^1\times Z_b$ and make the
obvious identification that $S(T(S^1\times Z_b))={\Bbb C}^2\otimes S(TZ_b)$.
$D^\prime =\{ D_b^\prime \}_{b\in B}$ then forms a family of Dirac operators
over $B$.

The connection $\tilde{\nabla}$ in (2.6) lifts naturally to
$\Gamma (S(T(S^1\times Z))\otimes E)$, which we still denote by $\tilde{\nabla}$.

Then by the fundamental result of Bismut [B], for any $t>0$, the following differential
form is a representative of ${\rm ch}(D^\prime)$,
$$\alpha _t(b)={\rm Tr}_s^{S^1\times Z} [\exp (-(\tilde{\nabla}+\sqrt{t}
D_{b,u}^\prime)^2)].\eqno (3.3)$$

Now for $\varepsilon >0$, set
$$D_{b,\varepsilon}^\prime = \left( \begin{array}{cc} 0 & \sqrt{\varepsilon} \frac{\pt}{\pt 
u} + D_{b,u} \\
-\sqrt{\varepsilon} \frac{\pt}{\pt u} + D_{b,u} & 0
\end{array} \right) \eqno (3.4)$$

Then for any $t>0,\ \varepsilon >0$, one clearly has that
$$\alpha _{t,\varepsilon}(b)={\rm Tr}_s[\exp (-(\tilde{\nabla}+\sqrt{t}
D_{b,u,\varepsilon}^\prime)^2)]\eqno (3.5)$$
is also a representative of ${\rm ch}(D^\prime)$.

Now taking the adiabatic limit $\varepsilon \rightarrow 0$.  By proceeding as
in [BF, Section 2f] and [B] in a very simple situation, one gets
$$\beta _t(b)=\lim_{\varepsilon \rightarrow 0} \alpha _{t,\varepsilon}(b)=
{1\over \sqrt{\pi}}\int_{S^1} {\rm Tr}^{Z,{\rm odd}}[\exp (-(du{\partial \over \partial u}
+\tilde{\nabla}+\sqrt{t} D_{b,u})^2)],\eqno (3.6)$$
which is also a representative of ${\rm ch}(D^\prime)$ for any $t>0$.

One verifies
$$(du{\partial \over \partial u}+\tilde{\nabla}+\sqrt{t}D_{b,u})^2=(
\tilde{\nabla}+\sqrt{t}D_{b,u})^2+\sqrt{t}du{\partial D_{b,u}\over \partial u}.
\eqno (3.7)$$
 From (3.6) and (3.7), one deduces that
$$\beta _t(b)={1\over \sqrt{\pi}}\int_{S^1} {\rm Tr}^{\rm odd}[\exp (-(
\tilde{\nabla}+\sqrt{t}D_{b,u})^2-\sqrt{t}du{\partial D_{b,u}\over \partial u})]$$
$$=-{1\over \sqrt{\pi}}\int_{S^1}{\partial \over \partial s}\{ {\rm Tr}^{\rm even}
[\exp (-(\tilde{\nabla}+\sqrt{t}D_{b,u})^2-s\sqrt{t}{\partial D_{b,u}\over
\partial u})]\}_{s=0}du.\eqno (3.8)$$

Now since we have a periodic family, and also that we have taken all $B_i=0$
$(i\geq 1$) in the superconnection formalism, we see that the $\hat{\eta}$-forms
at the end points (or the same point) are the same, so they cancel in (2.35).

Compare (3.8) with (2.22), and by (2.25), we get finally the following result.

$\ $

{\bf Theorem 3.2}. {\it The following identity holds in $H^*(B)$,
$${\rm ch}({\rm ind}\, D^\prime)={\rm ch}({\rm sf}\{ D\}).\eqno (3.9)$$}

{\bf Remark 3.3}. Theorem 3.2 generalizes the well-known result of 
Atiyah-Patodi-Singer [APS1] to higher dimensions.

$\ $

{\bf Remark 3.4}. It was pointed in  ([DZ3]) that (3.9) should still hold on the $K$-theoretic level.
That is, one should have
$${\rm ind}\, D^\prime={\rm sf}\{ D\}\ \ {\rm in}\ \ K(B).\eqno (3.10)$$
Fangbing Wu  [W] has since then proven this true.

$\ $

Now if $B$ is also spin,  for any $u\in S^1$, one can define the total
Dirac type operator $D_u^X$ on $X$. Then one has also a spectral flow
${\rm sf}\{ D^X\}$. The following result is a relation between this spectral flow 
and the higher spectral flow of the fibered Dirac operators.

$\ $

{\bf Corollary 3.5}. {\it The following identity holds,
$${\rm sf}\{ D_{u\in S^1}^X\}=\int_B\hat{A}(TB){\rm ch}({\rm sf}\{ D_{u\in S^1}^Z
\}),\eqno (3.11)$$
where now the Chern character is the usual Chern character, not the normalized one (i.e. here one has the factor $2\pi \sqrt{-1}$).}

$\ $

{\bf Remark 3.6}. The proof of Theorem 3.2 above is not the same as the one indicated in 
[DZ3]. To pursue the proof in [DZ3], one first reduces the problem to the
case of (standard) Dirac operators family and then uses the Bismut superconnection
formalism to evaluate the Chern character. Since we will encounter the Bismut
superconnection in the next section, so we think it is better to avoid it in this section.

\sect{ Toeplitz families and higher spectral flows: a heat kernel computation
for the index bundle of Toeplitz families}

In this section we introduce what we call Toeplitz families which extend the
usual concept of Toeplitz operators to fibration case. Just as the index of
Toeplitz operators can be expressed via spectral flow, we will intepret the index
bundle of a Toeplitz family via higher spectral flow. This in turn allows us to 
give an evaluation of the Chern character of the index bundle of Toeplitz
families via $\hat{\eta}$-forms. In particular, when the operators are of Dirac type,
we get an explicit local index computation of the Chern character, which can be viewd
as an odd analogue of the Bismut local index theorem [B] for families of Dirac
operators on even dimensional manifolds. Besides the Bismut superconnection, which 
should
be involved naturally here in our odd analogue, we will also use a conjugation of it.
See d) and e) for details.

This section is organized as follows. In a), we introduce what we call a Toeplitz
family and give the formula for the Chern character of its index bundle. In b), we establish 
a relationship between the index bundle of Toeplitz families and
a natually associated higher spectral flow. In c), we give a heat kernel
formula computing the Chern character of the index bundle of Toeplitz families.
In d), we introduce a conjugate form of the Bismut superconnection [B] and
finally in e), we give the local index evaluation of the Chern character of the
index bundle of Toeplitz families.

$\ $

{\bf a). Toeplitz families}

$\ $

We consider the same geometric objects as in Sections 1 and 2. That is, we have a fibration 
$Z\rightarrow M\stackrel{\pi}{\rightarrow} B$ with compact closed odd dimensional
fibers and compact base $B$. We assume the vertical bundle
$TZ\rightarrow M$ is spin and carries a 
fixed spin structure. $E\rightarrow M$ a complex  vector bundle with a
Hermitian metric $g^E$ and a Hermitian connection $\nabla ^E$. Also $TZ$ has
a metric $g^{TZ}$. So we have a canonically defined $B$-family of self-adjoint
(twisted) Dirac operators $D^E=\{ D_b^E\}_{b\in B}$.

Recall that we have also made the basic assumption that the associated Atiyah-Singer index
bundle
$${\rm ind}\, D^E=0\ \ {\rm in}\ \ K^1(B).\eqno (4.1)$$

Let $D=\{ D_b\}_{b\in B}$ be a smooth $B$-family of self-adjoint Dirac type operators (
that is, for any $b\in B$, the symbol of $D_b$ is the same as that of $D_b^E$).
Recall that by Remark 1.6, the $K^1$-index of $D$ also vanishes.

By the fundamental result of Melrose-Piazza [MP], which we already recalled in
Section 1, there is then a spectral section $P=\{ P_b\}_{b\in B}$ of the family
$D$.

Now let $g:M\rightarrow GL(N,{\Bbb C})$ be a smooth map. Then 
$g$ can be viewed as an automorphism of the trivial
complex vector bundle ${\Bbb C}^N\rightarrow M$ over $M$.
For any $b\in B$, $g$ induces an automorphism  of ${\Bbb C}^N\rightarrow Z_b$. Thus
for any $b\in B$, $g$ induces a bounded map $g_b$ from
$L^2(S(TZ_b)\otimes E_b\otimes {\Bbb C}^N)$ to itself by acting as an identity on
$L^2(S(TZ_b)\otimes E_b)$.

Also it is clear that $P_b$ induces a map on 
$P_bL^2(S(TZ_b)\otimes E_b\otimes {\Bbb C}^N)$ by acting as identity on sections of ${\Bbb C}^N$.
We still note this map by $P_b$.

$\ $

{\bf Definition 4.1}. The Toeplitz family $T_g=\{ T_{g,b}\}_{b\in B}$ is a family of
Toeplitz operators defined by
$$T_{g,b}=P_bg_bP_b:P_bL^2(S(TZ_b)\otimes E_b\otimes {\Bbb C}^N)\rightarrow
P_bL^2(S(TZ_b)\otimes E_b\otimes {\Bbb C}^N).\eqno (4.2)$$

$\ $

Clearly, $T_g$ is a continuous family. In fact one can even make it to be a smooth family
(cf. [MP]). Furthermore, since $g$ is invertible, for any $b\in B$, $T_{g,b}$ is a
Fredholm operator. Thus we get a family of Fredholm operators. By 
Atiyah-Singer [AS1] it then defines an index bundle
$${\rm ind}\, T_g\in K(B).\eqno (4.3)$$

Using a trick due to Baum and Douglas [BD], we now show that the Atiyah-Singer family 
index theorem
[AS1] provides a topological formula for this analytically defined ${\rm ind}\, T_g$.

In fact, recall that by [MP], $P_b$ and $g_b$ are zeroth order pseudodifferential
operators. So is $T_{g,b}$ for any $b\in B$.

Now for any $b\in B$, set
$$\tilde{T}_{g,b}=I-P_b+P_bg_bP_b:L^2(S(TZ_b)\otimes E_b\otimes {\Bbb 
C}^N)\rightarrow
L^2(S(TZ_b)\otimes E_b\otimes {\Bbb C}^N).\eqno (4.4)$$
Then by a simple calculation of symbol, we see that each $\tilde{T}_{g,b}$ is a
zeroth ordre elliptic pseudodifferential operator. Thus $\tilde{T}_g=\{ \tilde{T}
_{g,b}\} _{b\in B}$ is a continuous family of elliptic pseudodifferential
operators.  By  Atiyah-Singer [AS1], ${\rm ind}\, \tilde{T}_g$ can be identified
with a topological index bundle defined also by Atiyah-Singer [AS1]. Furthermore,
as we have  obviously
$${\rm ind}\,\tilde{T}_g={\rm ind}\, {T}_g\ \ {\rm in}\ \ K(B),\eqno (4.5)$$
we see that ${\rm ind}\, T_g$ is also equal to the above mentioned topological
index. In particular, its Chern character has a very natural form.

$\ $

{\bf Theorem 4.2}.  {\it The following identity holds,
$${\rm ch}({\rm ind}\, T_g)=-\int_{Z}\hat{A}(TZ){\rm ch}(E){\rm ch}(g)\ \ 
{\rm in}\ \ H^*(B),\eqno (4.6)$$
where ${\rm ch}(g)$ is the odd Chern character associated to $g$ (cf. [AH]).}

$\ $

In the rest of this section we will give a purely analytic proof of (4.6) via 
heat kernel method.

$\ $

{\bf b). Toeplitz families and higher spectral flows}

$\ $

We make the same assumption and use the same notation as in a). It is clear that 
each $D_b$ can also be extended as an operator from $\Gamma (S(TZ_b)\otimes 
E_b\otimes {\Bbb C}^N|_{Z_b})$
to itself by acting as identity on ${\Bbb C}^N$. 

Also, as $GL(N,{\Bbb C})$ is homotopically equivalent to $U(N)$, and that 
${\rm ind}\, T_g$ is clearly a homotopic invariant of $g\in GL(N,{\Bbb C})$, 
without loss of generality in the computation of ${\rm ind}\, T_g$, from now on
we will assume $g\in U(N)$.

Set for any $b\in B$,
$$D_{g,b}=g_bD_bg_b^{-1}.\eqno (4.7)$$
Since $g\in U(N)$, one verifies easily that $D_g=gDg^{-1}=\{ D_{g,b}\}_{b\in B}$ is also
a $B$-family of self-adjoint Dirac type operators.

For any $0\leq t\leq 1$, set
$$D_b(t)=(1-t)D_b+tD_{g,b}.\eqno (4.8)$$
Then we get a natural curve $D(t)=\{ D_b(t)\}_{b\in B},\ t\in [0,1]$ of $B$-families
of self-adjoint Dirac type operators.

Let $P$ be a spectral section of $D$. Then clearly $gPg^{-1}=\{ g_bP_bg_b^{-1}\}
_{b\in B}$ is a spectral section of $gDg^{-1}$. And we have a naturally defined higher 
spectral flow
$${\rm sf}\{ (D,P),(gDg^{-1},gPg^{-1})\}\in K(B)\eqno (4.9)$$
associated to the curve (4.8) (See Section 1).

$\ $

{\bf Proposition 4.3}. {\it i), ${\rm sf}\{ (D,P),(gDg^{-1},gPg^{-1})\}$ does not
depend on $P$; ii), ${\rm sf}\{ (D,P),(gDg^{-1},gPg^{-1})\}$ depends only on the
symbol of $D$.}

{\it Proof}. i). Let $Q$ be another spectral section of $D$. Then one verifies easily
that
$${\rm sf}\{ (D,P),(gDg^{-1},gPg^{-1})\} 
-{\rm sf}\{ (D,Q),(gDg^{-1},gQg^{-1})\}$$ 
$$=[gPg^{-1}-gQg^{-1}]-[P-Q]\in K(B).\eqno (4.10)$$
Now by Kuiper's theorem [K], as a map into the unitary group of a Hilbert space, $g$ is
homotopic to $Id$ via unitary operators. Thus one has, by the homotopy 
invariance of the $K$-group
elements, that
$$[gPg^{-1}-gQg^{-1}]=g[P-Q]g^{-1}=[P-Q]\ {\rm in}\ K(B).\eqno (4.11)$$
The conclusion of i) follows from (4.10) and (4.11).

ii). Let $D^\prime$ be another $B$-family of self-adjoint Dirac type operators
with the same symbol as $D$. Take any spectral section $R$ of $D^\prime$. Since
$D$, $D^\prime$ has the same symbol, a spectral section of $D$ is  a generalized
spectral section of $D^\prime$ and vise versa. Furthermore, since for each $t\in [0,1]$,
$D(t)$ defined in (4.8) has the same symbol as $D$, one deduces that
$${\rm sf}\{ (D,P),(gDg^{-1},gPg^{-1})\}-{\rm sf}\{ (D^\prime, R),(gD^\prime g^{-
1},
gRg^{-1})\}$$
$$=[gPg^{-1}-gRg^{-1}]-[P-R]\ \ {\rm in}\ \ K(B),\eqno (4.12)$$
which vanishes by the same reason as in the proof of i). \qed

$\ $

Thus, what we get is an intrinsic invariant depending only on $g$ and the symbol
of $D^E$. We will denote it by ${\rm sf}\{ D^E,gD^Eg^{-1}\}$ in what follows.

We now state the main result of this subsection, which can be seen as an extension
to fibration case of the classical result of Booss-Wojciechowski (cf. [BW, Theorem 
17.17]).

$\ $

{\bf Theorem 4.4}. {\it The following identity holds in $K(B)$,
$${\rm ind}\, T_g={\rm sf}\{ D^E,gD^Eg^{-1}\}.\eqno (4.13)$$}

{\it Proof}. The proof of (4.13) is already implicitely contained in the proof of
ii) of Proposition 4.3. The main observation is that for any $t\in [0,1]$, the
symbol of $D(t)$ is the same as that of $D^E$. Thus $P$ is a generalized spectral
section in the sense of Definition 1.8  for all $D(t),\ 0\leq t\leq 1$.
By the (modified) definition of higher spectral flow, one then gets
$${\rm sf}\{ D^E,gD^Eg^{-1}\} =[gPg^{-1}-P]={\rm ind}\, T_g\ \in\ K(B).
\eqno (4.14)$$

The proof of Theorem 4.4 is completed. \qed

$\ $

{\bf c). A heat kernel formula for the Chern character of the index bundle of
Toeplitz families}

$\ $

We combine Theorem 4.4 and Theorem 2.12 to give a heat kernel formula for 
${\rm ind}\, T_g$. This amounts first to introduce as in Section 2 an
orthogonal splitting
$$TM=TZ\oplus T^HM,$$
$$g^{TM}=g^{TZ}\oplus \pi ^*g^{TB}.\eqno (4.15)$$

We now follow the same strategy as in Section 2 by taking $D_1=gDg^{-1}$, and
consider the higher spectral flow ${\rm sf}\{ (D,P),(gDg^{-1},gPg^{-1})\}$.

For our special situation, we make the assumption that in the notation of (2.21),
one has
$$B(t,1)=gB(t,0)g^{-1}.\eqno (4.16)$$
This clearly can always be achieved.

By (4.16), one gets immediately
$$\hat{\eta} (B_1,gPg^{-1})=\hat{\eta} (B_0,P).\eqno (4.17)$$
Thus by Theorem 2.12, one gets

$\ $

{\bf Proposition 4.5}. {\it The following identity holds in $\Omega (B)/d\Omega
(B)$,
$${\rm ch}({\rm ind}\, T_g)=-{1\over \sqrt{\pi}}\int_0^1 a_0(u)du,\eqno (4.18)$$
where $a_0(u)$ is the asymptotic term in (2.22) subject to the condition
(4.16) of the superconnection $B(t,u)$.}

{\it Proof}. By (4.17) and Theorem 2.12, one has
$${\rm ch}({\rm sf}\{ D^E,gD^Eg^{-1}\})=-{1\over \sqrt{\pi}}\int_0^1 a_0(u)du.
\eqno (4.19)$$

(4.18) follows from (4.19) and Theorem 4.4. \qed

$\ $

{\bf d). A conjugate form of the Bismut superconnection}

$\ $

We use the same notation as in Section 2c). Also, following Bismut [B], let $S$
be the tensor defined by
$$\nabla ^{TM}=\nabla ^{TZ}+\pi ^*\nabla ^{TB}+S.\eqno (4.20)$$

Clearly, the connection $\tilde{\nabla}$ defined in (2.6) extends canonically
as a connection on $H_{\infty}\otimes {\Bbb C}^N$.

Let $e_1,\ldots ,e_n$ be an orthonormal base of $TZ$. Set
$$k=-{1\over 2}\sum_1^n S(e_i)e_i,\eqno (4.21)$$
and let $\tilde{\nabla}^u$ be the connection on $H_{\infty}$ defined by
$$\tilde{\nabla}_Y^u=\tilde{\nabla}_Y+\langle k,Y^H\rangle,\ Y\in TB.\eqno (4.22)$$
Then $\tilde{\nabla}^u$ is a unitary connection on $H_{\infty}\otimes {\Bbb C}^N$,
and does not depend on $g^{TB}$.

Also let $T$ be the torsion of $\nabla ^{TZ}+\pi ^*\nabla ^{TB}$. Set
$$c(T)=-\sum_{\alpha <\beta} \langle T(f_{\alpha},f_\beta),e_i\rangle c(e_i)
dy^\alpha dy^\beta,\eqno (4.23)$$
where $\{ f_\alpha \}$ is an orthonormal base of $TB$ and $\{ dy^\alpha \}$ the
dual base of $T^*B$ (We identify $TB$ with $\pi ^*TB$).

$\ $

{\bf Definition 4.6} ([B]). The Bismut superconnection $B(0)$ is defined by
$$B(0)=\tilde{\nabla} ^u+D^E-{c(T)\over 4}.\eqno (4.24)$$
The rescaled Bismut superconnection is then given by
$$B_t(0)=\tilde{\nabla} ^u+\sqrt{t}D^E-{c(T)\over 4\sqrt{t}}.\eqno (4.25)$$

Recall that $g$ is unitary as in b).

$\ $

{\bf Definition 4.7}. The $g$-conjugate Bismut superconnection $B(1)$ is defined by
$$B(1)=gB(0)g^{-1}.\eqno (4.26)$$
The rescaled $g$-conjugate Bismut superconnection $B_t(1)$ is then given by
$$B_t(1)=gB_t(0)g^{-1}.\eqno (4.27)$$

For any $u\in [0,1]$, set
$$B_t(u)=(1-u)B_t(0)+uB_t(1).\eqno (4.28)$$
We will use $B_t(u)$ to calculate the asymptotic term $a_0(u)$ in (4.18). The
result will be stated in Theorem 4.8, which implies an explicit evaluation of
${\rm ch}({\rm ind}\, T_g)$.

$\ $

{\bf e). A local formula for the Chern character of the index bundle of 
Toeplitz families}

$\ $

We will use the superconnections in d) to calculate the asymptotic term $a_0(u)$
in (4.18).

Recall that $a_0(u)$ is defined in (2.22) by an asymptotic formula. In using
notation in d),
$${\partial \over \partial b}\{ {\rm Tr}^{\rm even}[\exp (-B^2_t(u)-b
{\partial B_t(u)\over \partial u})]\} _{b=0}$$
$$={a_{-k}\over t^k}+\cdots +{a_{-1}\over t}+a_0+o(1),\ t\rightarrow 0.\eqno (4.29)$$

The following is the main result of this section.

$\ $

{\bf Theorem 4.8}. {\it The following identity holds,
$$\lim_{t\rightarrow 0} {\rm Tr}^{\rm even}[{\partial B_t(u)\over \partial u}
\exp (-B^2_t(u))]$$
$$=-({1\over 2\pi \sqrt{-1}})^{\dim Z+1\over 2}\sqrt{\pi}\int_Z\hat{A}(R^{TZ})
{\rm Tr}[\exp (-R^E)]{\rm Tr}[g^{-1}dg\exp ((u-u^2)(g^{-1}dg)^2)],\eqno (4.30)$$
where $R^{TZ}$ is the curvature of $\nabla ^{TZ}$ and $R^E$ is the curvature of
$\nabla ^E$.}

{\it Proof}. By (4.27), one has
$$B_t(1)=gB_t(0)g^{-1}=B_t(0)+g[B_t(0),g^{-1}].\eqno (4.31)$$

Thus for any $u\geq 0$,
$$B_t(u)=B_t(0)+ug[B_t(0),g^{-1}]\eqno (4.32)$$
and
$${\partial B_t(u)\over \partial u}=g[B_t(0),g^{-1}].\eqno (4.33)$$

Also one verifies easily that
$$g[B_t(0),g^{-1}]=g[\tilde{\nabla} ^u+\sqrt{t}D-{c(T)\over 4\sqrt{t}},g^{-1}]$$
$$=g[\tilde{\nabla}^u+\sqrt{t}D,g^{-1}]=g\sum_\alpha dy^\alpha ({\partial \over
\partial y^\alpha}g^{-1})+g\sum_i\sqrt{t}c(e_i)({\partial \over \partial e_i}
g^{-1})$$
$$=-\sum_\alpha dy^\alpha ({\partial \over \partial y^\alpha}g)g^{-1}-\sum_i\sqrt{t}
c(e_i)({\partial \over \partial e_i}g)g^{-1}.\eqno (4.34)$$

Also from (4.32), one deduces that
$$B_t^2(u)=B_t^2(0)+u[B_t(0),g[B_t(0),g^{-1}]]+u^2g[B_t(0),g^{-1}]g[B_t(0),g^{-
1}]$$
$$=B_t^2(0)+u[B_t(0),g][B_t(0),g^{-1}]+ug[B_t(0),[B_t(0),g^{-
1}]]+u^2(g[B_t(0),g^{-1}])^2.
\eqno (4.35)$$

Now one verifies by (4.34) that
$$[B_t(0),g]=\sum_\alpha dy^\alpha {\partial g\over \partial y^\alpha}+\sqrt{t}
\sum_i c(e_i){\partial g\over \partial e_i}$$
$$=-g[B_t(0),g^{-1}]g.\eqno (4.36)$$

If we make the formal change that $\sqrt{t}c(e_i)\rightarrow e_i\in T^*Z$, then
one gets
$$[B_t(0),g^{-1}]\rightarrow dg^{-1}\eqno (4.37)$$
and
$$[B_t(0),[B_t(0),g^{-1}]]\rightarrow d^2g^{-1}=0.\eqno (4.38)$$

 From (4.33), (4.35) to (4.38), and by proceeding as in [B], [BF], and 
[BGV, Chap. 10] the by now
standard local index techniques, which can be adapted easily here in odd 
dimensions, one gets that as $t\rightarrow 0$, one has
$$\lim_{t\rightarrow 0} {\rm Tr}^{\rm even}[{\partial B_t(u)\over \partial u}\exp 
(-B_t^2(u))]$$
$$=({1\over 2\pi \sqrt{-1}})^{\dim Z+1\over 2}\sqrt{\pi}\int_Z\hat{A}(R^{TZ})
{\rm Tr}[\exp (-R^E)]{\rm Tr}[gdg^{-1}\exp ((u-u^2)(gdg^{-1})^2)].\eqno (4.39)$$

Now clearly 
$$gdg^{-1}=-(dg)g^{-1}=-g(g^{-1}dg)g^{-1}.\eqno (4.40)$$

(4.30) follows from (4.39) and (4.40). \qed

$\ $

 From (4.30), one finds
$${1\over \sqrt{\pi}}\int_0^1\lim_{t\rightarrow 0} {\rm Tr}^{\rm even}[
{\partial B_t(u)\over \partial u}\exp (-B_t^2(u))]du$$
$$=-({1\over 2\pi \sqrt{-1}})^{\dim Z+1\over 2}\int_Z\hat{A}(R^{TZ}){\rm Tr}[\exp
(-R^E)]\sum_{n=0}^\infty {n!\over (2n+1)!}{\rm Tr}[(g^{-1}dg)^{2n+1}].
\eqno (4.41)$$

We can now state the main result of this section, which is a local form of
Theorem 4.2.

$\ $

{\bf Theorem 4.9}. {\it The following differential form on $B$,
$$-({1\over 2\pi \sqrt{-1}})^{\dim Z+1\over 2}\int_Z\hat{A}(R^{TZ}){\rm Tr}[\exp
(-R^E)]\sum_{n=0}^\infty {n!\over (2n+1)!}{\rm Tr}[(g^{-1}dg)^{2n+1}],
\eqno (4.42)$$
is closed and represents ${\rm ch}({\rm ind}\, T_g)$ in $H^*(B)$.  \qed}

$\ $

{\bf Remark 4.10}. In the definition of Toeplitz families, one can in fact replace
$g$ by any automorphism of an arbitary vector bundle over $M$. From the 
topological point of view, however,  the trivial bundle is sufficiently general,  for $H^{\rm odd}(M)$
is generated by the odd Chern characters of automorphisms of trivial vector
bundles (cf. [AH]).

\sect{Higher spectral flows and family indices for manifolds with boundary}

In this section we use the idea of higher spectral flow to study the index theory
for families of manifolds with boundary. Among the two main results of this
section, one generalizes an earlier result in [DZ1, Theorem 1.1] to family case. It 
characterizes
the changes of family indices under continuous deformations of Dirac type operators.
The other one extends the Melrose-Piazza index theorem [MP] to the general
superconnection case, not necessary restricted to the Bismut superconnections. Of
course then, the local index density can not be directly identified with
characteristic forms.

This section is organized as follows. In a), we define what we call a $B$-family of
Dirac type operators on (even) dimensional manifolds with boundary. And following
[MP], for any spectral section on the boundary family, we define an index bundle
of the $B$-family of Dirac type operators. In b), we prove our result on the
change of index bundles under various deformations. In c), we prove the index formula mentioned above.

$\ $

{\bf a). A family of Dirac type operators on even dimensional manifolds with
boundary}

$\ $

Let $Y\rightarrow M\stackrel{\pi_Y}{\rightarrow} B$ be a smooth fibration with the typical fibre
an even dimensional compact manifold $Y$ with boundary $Z$. We assume $B$ is
compact. 
Again, for simplicity we make the assumption that the vertical bundle
$TY$ is spin and carries a fixed
spin structure. Then the  boundary fibration $Z\rightarrow \partial M\stackrel{\pi _Z}{\rightarrow} B$
verifies the properties of previous sections.

Let $g^{TY}$ be a metric on $TY$. Let $g^{TZ}$ be the restriction of $g^{TY}$ on
$TZ$. We make the assumption that there is a neighborhood $[-1,0]\times \partial M$
of $\partial M$ in $M$ such that for any $b\in B$, on 
$\pi _Y^{-1}(b)\cap ({[-1,0]\times
\partial M})=[-1,0]\times Z_b$, $g^{TY}$ takes the form
$$g^{TY}|_{[-1,0]\times Z_b}=dt^2\oplus g^{TZ_b}.\eqno (5.1)$$

We denote by $\pi _b:[-1,0]\times Z_b\rightarrow Z_b$ the projection onto the
second factor.

Let $S(TY)$ be the spinor bundle of $(TY,g^{TY})$. Since $\dim Y$ is even, we have
a canonical splitting of $S(TY)$ to positive/negative spinor bundles:
$S(TY)=S_+(TY)\oplus S_-(TY)$.

Let $E$ be a complex vector bundle over $M$. Let $g^E$ be a metric on $E$ such that
$$g^E|_{[-1,0]\times Z}=\pi ^*g^E|_Z.\eqno (5.2)$$
Let $\nabla ^E$ be a Hermitian connection on $E$ such that 
$$\nabla ^E|_{[-1,0]\times Z}=\pi ^*\nabla ^E|_Z.\eqno (5.3)$$
The exsitence of such a pair of $(g^E,\nabla ^E)$ is clear.

Let $S(TZ)$ be the spinor bundle of $(TZ,g^{TZ})$. We will still use the notation
in previous sections for $Z$ here. Thus there are canonically defined (twisted)
Dirac opeartors $D_b^{Z,E}$ associated to $g^{TZ_b},\ g^{E|_{Z_b}}$ and $\nabla
^{E|_{Z_b}}$.

Since $Z\rightarrow \partial M\rightarrow B$ is a boundary family, it is known
([Sh], see also [N] for a new proof) that the $K^1(B)$ index of $D^{Z,E}=\{ D_b^{Z,E}\}_{b\in B}$ vanishes.

Let $P$ be a spectral section of $D^{Z,E}$, the existence
of which is guaranteed by [MP].

Now for any $b\in B$, there is also a (twisted) Dirac operator
$D_b^{Y,E}:\Gamma (S_+(TY_b)\otimes E|_{Y_b})\rightarrow \Gamma (S_-
(TY_b)\otimes
E|_{Y_b})$ canonically associated to $g^{TY_b},\ g^{E|_{Y_b}}$ and 
$\nabla ^{E|_{Y_b}}$. Clearly on $U_b=[-1,0]\times Z_b$, $D_b^{Y,E}$ takes the
form
$$D_b^{Y,E}=c({\partial \over \partial t})({\partial \over \partial t}+D_b^{Z,E}).\eqno 
(5.4)$$

$\ $

{\bf Definition 5.1}. By a $B$-family of Dirac type operators on $Y\rightarrow 
M\rightarrow B$,
we will mean a $B$-family of first order differential operators
$D^Y=\{ D_b^Y\} _{b\in B}$ with each $D_b^Y:\Gamma (S_+(TY_b)\otimes E|_{Y_b})
\rightarrow \Gamma (S_-(TY_b)\otimes E|_{Y_b})$ has the same symbol as that of
$D_b^{Y,E}$ and such that on $U_b=[-1,0]\times Z_b$, $D_b^Y$ can be written in
the form
$$D_b^Y=c({\partial \over \partial t})({\partial \over \partial t}+D_b^Z),
\eqno (5.5)$$
where $D^Z=\{ D_b^Z\}_{b\in B}$ is a $B$-family of (self-adjoint) Dirac type operators
of $Z\rightarrow \partial M\rightarrow B$.

$\ $

Recall that $P$ is a spectral section of $D^{Z,E}$. It is then a generalized spectral section
of $D^Z$ in the sense of Definition  . For any $b\in B$, we then have a generalized
Atiyah-Patodi-Singer elliptic boundary problem $(D_b^Y,P_b)$. And the family
$(D^Y,P)=\{ (D_b^Y,P_b)\} _{b\in B}$ gives a $B$-family of Fredholm operators.
According to Atiyah and Singer [AS1], it then determines an index bundle
$${\rm ind}\, (D^Y,P)\in K(B).\eqno (5.6)$$

$\ $

{\bf b). A variation formula for index bundles}

$\ $

Now consider a curve of $B$-families of metrics and connections
$(g^{TY,u},g^{E,u},\nabla ^{E,u})$, $u\in [0,1]$ verifying the conditions
(5.1) to (5.3) for each $u\in [0,1]$.

Let $\tilde{D}=\{ D^Y(u)\}_{u\in [0,1]}$ be a smooth curve of $B$-family of 
Dirac type operators such that for any $u$, $D^Y(u)$ has the same symbol as that of
$D^{Y,E}(u)$ corresponding to $(g^{TY,u},g^{E,u},\nabla ^{E,u})$.

Let $P_0,\ P_1$ be two generalized spectral sections of $D^Z(0),\ D^Z(1)$ 
respectively. Then as in a), one has two naturally defined index bundles
$${\rm ind}\, (D_Y(0),P_0)\in K(B),$$
$${\rm ind}\, (D^Y(1),P_1)\in K(B).\eqno (5.7)$$

On the other hand, $(D^Z(0),P_0),\ (D^Z(1),P_1)$ and the curve
$\{ D^Z(u)\} _{u\in [0,1]}$ determines a higher spectral flow by Section 1:
$${\rm sf}\{ (D^Z(0),P_0),(D^Z(1),P_1)\} \in K(B).\eqno (5.8)$$

The first main result of this section, which generalizes 
[DZ1, Theorem 1.1], can be stated as follows.

$\ $

{\bf Theorem 5.2}. {\it The following identity holds in $K(B)$,
$${\rm ind}\, (D^Z(0),P_0)-{\rm ind}\, (D^Z(1),P_1)={\rm sf}\{ (D^Z(0),P_0),
(D^Z(1),P_1)\} .\eqno (5.9)$$}

{\it Proof}. Let $\tilde{Q}$ be a generalized spectral section of the total 
family $\tilde{D}^Z=\{ D^Z_b(u)\}$. Then we get a continuous
family of elliptic boundary problems $(D^Y(u),\tilde{Q}(u))_{u\in [0,1]}$.

Thus by the homotopy invariance of the index bundles, one gets
$${\rm ind}\, (D^Y(0),\tilde{Q}(0))={\rm ind}\, (D^Y(1),\tilde{Q}(1))\ {\rm in}\ 
K(B).\eqno (5.10)$$

Now by the relative index theorem in [DZ2], which extends the relative index 
theorem of Melrose-Piazza [MP] to generalized spectral sections, one has
$${\rm ind}\, (D^Y(0),\tilde{Q}(0))={\rm ind}\, (D^Y(0),P_0)+[P_0-\tilde{Q}(0)]\ 
{\rm in}\ K(B)\eqno (5.11)$$
and
$${\rm ind}\, (D^Y(1),\tilde{Q}(1))={\rm ind}\, (D^Y(1),P_1)+[P_1-\tilde{Q}(1)]\ 
{\rm in}\ K(B).\eqno (5.12)$$

 From (5.10) to (5.12), one gets
$${\rm ind}\, (D^Y(0),P_0)-{\rm ind}\, (D^Y(1),P_1)=[P_1-\tilde{Q}(1)]-[P_0-
\tilde{Q}(0)]$$
$$={\rm sf}\{ (D^Z(0),P_0),(D^Z(1),P_1)\}\ {\rm in}\ K(B).\eqno (5.13)$$

The proof of Theorem 5.2  is completed. \qed

$\ $

{\bf c). A heat kernel asymptotic formula for index bundles}

$\ $

In this subsection, we give a heat kernel asymptotic formula for the index bundle
(5.6) in terms of a local index density in the interior and a $\hat{\eta}$-form
on the boundary. Since the $\hat{\eta}$ form is only well-defined with repect to
spectral sections, not the generalized ones, we will restric in this subsection
only to consider the former ones.

We first introduce the superconnection formalism.

For the boundary fibration $Z\rightarrow \partial M\rightarrow B$, we choose
a splitting as in Section 2:
$$T\partial M=TZ\oplus \pi ^*_ZTB,$$
$$g^{T\partial M}=g^{TZ}\oplus \pi ^*_Zg^{TB}.\eqno (5.14)$$

This splitting induces a trivial product splitting on $U=[-2,0]\times \partial M$:
$$TU=T([-2,0])\oplus T\partial M,$$
$$g^{TU}=g^{T([-2,0])}\oplus g^{T\partial M}.\eqno (5.15)$$

We assume $(TM,g^{TM}$) also carries a splitting
$$TM=TY\oplus \pi ^*_YTB,$$
$$g^{TM}=g^{TY}\oplus \pi^*_Yg^{TB}\eqno (5.16)$$
with
$$(TM,g^{TM})|_U=(TU,g^{TU}).\eqno (5.17)$$

Let $\tilde{M}$ be the double of $M$ along $\partial M$ via the product structure near
$\partial M$. Then $\tilde{M}$ is a fibration $\tilde{Y}\rightarrow \tilde{M}\stackrel{\tilde{\pi}}{\rightarrow}
B$
such that for any $b\in B$, $\tilde{\pi}^{-1}(b)=\tilde{Y}_b$ is the double of $Y_b$
via the product structure near $Z_b=\partial Y_b$. The splittings (5.16), (5.17)
determine a splitting of $(T\tilde{M},g^{T\tilde{M}})$.

Let $D^Y=\{ D^Y_b\}_{b\in B}$ be a $B$-family of Dirac type operators for
$\pi _Y$. Then for any $b\in B$, $D^Y_b$ extends canonically to an invertible
Dirac type operator on $\tilde{Y}_b$ (cf. [BW]). Thus we get a $B$-family of
invertible Dirac type operators $D^{\tilde{Y}}=\{ D_b^{\tilde{Y}}\}_{b\in B}$.

Now the splitting (5.16) determines a connection $\tilde{\nabla}$ of the infinite
dimensional bundle
$$H_\infty (\tilde{Y},E)=\Gamma (S(T\tilde{Y})\otimes E)\eqno (5.18)$$
over $B$ by (2.6), which preserves the natural splitting
$$\Gamma (S(T\tilde{Y})\otimes E)=H_{\infty, +}\oplus H_{\infty, -}$$
$$=\Gamma (S_+(T\tilde{Y})\otimes E)\oplus \Gamma (S_-(T\tilde{Y})\otimes E).
\eqno (5.19)$$

Set $\tilde{U}=[-2,2]\times \partial M$, the product neighborhood of $\partial M$
in $\tilde{M}$ induced canonically by $U=[-2,0]\times \partial M\subset M$. The
metric splitting and connection splitting clearly holds on $\tilde{U}$. In
particular, $D^{\tilde{Y}}$ is of product structure near $\tilde{U}$. 

Now let $B_i^Y\in \Omega ^i(B)\hat{\otimes }\Gamma ({\rm cl}(TY)\otimes {\rm End }E),\ 
i\geq 1$ be odd with respect to the natural ${\bf Z}_2$ grading of 
$\Lambda (T^*B)\otimes {\rm cl}(TY)$, and that $B_i^Y|U$ does not depend on
$u\in [-2,0]$. Then $B_i^Y$ extends to $\tilde{Y}$  as
$B_i^{\tilde{Y}}$, the restriction of which to $\tilde{U}$ does not depend on
$u\in [-2,2]$.

For any $t>0$, set
$$B_t^{\tilde{Y}}=\tilde{\nabla} +\sqrt{t} D^{\tilde{Y}}+\sum_{i\geq 1}
t^{1-i\over 2}B_i^{\tilde{Y}}.\eqno (5.20)$$
Then $B_t^{\tilde{Y}}$ is a superconnection on $H_\infty (\tilde{Y},E)$ in the
sense of Quillen [Q] and Bismut [B].

One verifies  easily that as $t\rightarrow 0$, one has the uniform asymptotic
expansion on $\tilde{Y}$,
$${\rm Tr}_s[\exp (-B_t^{\tilde{Y},2})(y,y)]=\sum_{i\geq -{\dim Y\over 2}-
[{\dim B\over 2}]}^0 a_it^i+o(1),\ y\in \tilde{Y},\eqno (5.21)$$
where $\exp (-B_t^{\tilde{Y},2})(x,y),\ x,\, y\in \tilde{Y}$ is the $C^\infty$
kernel of $\exp (-B_t^{\tilde{Y},2})$ along the fibre $\tilde{Y}$.

On the other hand, each $B_i^Y$ also induces a natural element $B_i^Z$ in
$\Omega ^i(B)\hat{\otimes} \Gamma ({\rm cl}(TZ)\otimes E|_{\partial M})$. Also 
$\tilde{\nabla}$ induces a connection $\tilde{\nabla}^Z$ on
$$H_\infty (Z,E|_Z)=\Gamma (S(TZ)\otimes E|_{\partial M})\eqno (5.22)$$
and one has the induce superconnection for any $t>0$
$$B_t^Z=\tilde{\nabla}^Z+\sqrt{t}D^Z+\sum_{i\geq 1}t^{1-i\over 2}B_i^Z.\eqno
(5.23)$$

Let $P$ be a spectral section of $D^Z$ in the sense of [MP]. Then by Section 2,
one has a well-defined $\hat{\eta}$ form (class)
$$\hat{\eta}(D^Z,P,B^Z)\in \Omega (B)/d\Omega (B).\eqno (5.24)$$

We can now state the second main result of this section, which extends the family index
theorems in [BC2], [BC3] and [MP].

$\ $

{\bf Theorem 5.3}. {\it The following identity holds in $H^*(B)$,
$${\rm ch}({\rm ind}\, (D^Y,P))=\int_Y a_0-\hat{\eta} (D^Z,P,B^Z),\eqno (5.25)$$
where $a_0$ is given by (5.21).}

{\it Proof}. Let $D^{Y,E}$ be the standard $B$-family of (twisted) Dirac operators
associated to $(g^{TY},g^E,\nabla ^E)$. Then
$$D^Y(s)=(1-s)D^Y+sD^{Y,E},\ \ 0\leq s\leq 1,\eqno (5.26)$$
forms a curve of $B$-families of Dirac type operators over $B$.

Let $Q$ be a spectral section of $D^{Z,E}$, which is the boundary family
induced by $D^{Y,E}$. Then by Theorem 5.2, one has
$${\rm ind}\, (D^Y,P)-{\rm ind}\, (D^{Y,E},Q)={\rm sf}\{ (D^Z,P),(D^{Z,E},Q)\}.
\eqno (5.27)$$

Now let $B^{Y,E}$ be the Bismut superconnection [B] associated to $D^{Y,E}$ and
the splitting (5.16). Let $B^{Z,E}$ be the induced Bismut superconnection on
the boundary family. Then it is a result of Melrose-Piazza [MP],
which generalizes the earlier results of Bismut-Cheeger [BC2], [BC3],  that (5.25)
holds for $(B^{Y,E},B^{Z,E})$. That is, one has the formula
$${\rm ch}({\rm ind}\, (D^{Y,E},Q))=\int_Y a_0^E-\hat{\eta}(D^{Z,E},Q,B^{Z,E}),\ 
{\rm in}\ H^*(B).\eqno (5.28)$$

Now for any $s\in [0,1]$, $t>0$,
$$B_t^Y(s)=(1-s)B_t^Y+sB_t^{Y,E}\eqno (5.29)$$
is a (rescaled) superconnection associated to the $B$-family of Dirac type 
operators $D^Y(s)$. Denote its extension to $\tilde{M}$ by $B_t^{\tilde{Y}}(s)$.

For any $t>0,\ s\in [0,1],\ y,\, y^\prime \in \tilde{\pi} ^{-1}(b)$, $b\in B$,
let $P_{t,s}(y,y^\prime )=\exp (-I_{t,s})(y,y^\prime )$ be the smooth kernel
associated to $I_{t,s}=(B_t^{\tilde{Y}}(s))^2$. By using the Duhamel principle
and by proceeding as in [B, Section 2], one deduces easily that
$${\partial \over \partial s}{\rm Tr}_s[P_{t,s}(y,y)]=-{\rm Tr}_s[
{\partial I_{t,s}\over \partial s}P_{t,s}(y,y^\prime )]_{y=y^\prime}$$
$$=-{\rm Tr}_s[\, [B_t^{\tilde{Y}}(s),{\partial B_t^{\tilde{Y}}(s)\over \partial s}]\exp
(-B_t^{\tilde{Y}}(s)^2)(y,y^\prime )]_{y=y^\prime}$$
$$=-{\rm Tr}_s[\, [B_t^{\tilde{Y}}(s),{\partial B_t^{\tilde{Y}}(s)\over
\partial s}\exp (-B_t^{\tilde{Y}}(s)^2)](y,y^\prime )]_{y=y^\prime}$$
$$=-d{\rm Tr}_s[{\partial B_t^{\tilde{Y}}(s)\over \partial s}\exp (-B_t^{\tilde{Y}}(s)^2)
(y,y^\prime )]_{y=y^\prime}$$
$$-\sqrt{t}{\rm Tr}_s[\, [D^{\tilde{Y},E},{\partial B_t^{\tilde{Y}}(s)\over \partial s}
\exp (-B_t^{\tilde{Y}}(s)^2)](y,y^\prime )]_{y=y^\prime}.\eqno (5.30)$$

Thus one has
$${\partial \over \partial s}{\rm Tr}_s[\exp (-B_t^{\tilde{Y}}(s)^2)(y,y )]$$
$$=-\sqrt{t}{\rm Tr}_s[\, [D^{\tilde{Y},E},{\partial B_t^{\tilde{Y}}\over \partial s}
\exp (-B_t^{\tilde{Y}}(s)^2)](y,y^\prime )]_{y=y^\prime}\ {\rm mod}\ 
(\tilde{\pi}d^B)\Omega (\tilde{Y}).\eqno (5.31)$$

Integrating (5.31) from $s=0$ to $s=1$, one gets
$${\rm Tr}_s[\exp (-B_t^{\tilde{Y},2})(y,y)]-
{\rm Tr}_s[\exp (-B_t^{\tilde{Y},E,2})(y,y)]$$
$$=\sqrt{t}\int_0^1{\rm Tr}_s[\, [D^{\tilde{Y},E},{\partial B_t^{\tilde{Y}}(s)\over
\partial s}\exp (-B_t^{\tilde{Y}}(s)^2)](y,y^\prime )]_{y=y^\prime}ds\ {\rm mod}\ 
(\tilde{\pi}d^B)\Omega (\tilde{Y}).\eqno (5.32)$$

Now we let $t\rightarrow 0$ in (5.32). By procedding as in [B] and by 
comparing the coefficients of the resulting asymptotics, we get
$$a_0-a_0^E=d^{\tilde{Y}}\int_0^1 b_0(s)ds\ {\rm mod}\ (\tilde{\pi}d^B)\Omega
(\tilde{Y}),
\eqno (5.33)$$
where $b_0(s)$ is the constant term in the asymptotics as $t\rightarrow 0$ of
$${\rm Tr}_s[{\partial B_t^{\tilde{Y}}(s)\over \partial s}\exp (-
B_t^{\tilde{Y}}(s)^2)(y,y^\prime )]_{y=y^\prime}.\eqno (5.34)$$

Also, by the product structure of metrics and connections near the boundary
fibration $\partial M$, we see that if $c_0$ is the constant term of the
asymptotics of
$${\rm Tr}^{\rm even}[{\partial B_t^Z(s)\over \partial s}\exp (-
B_t^Z(s)^2)(z,z^\prime )]_{z=z^\prime},\ z,\, z^\prime \in Z,\eqno (5.35)$$
where $B_t^Z(s)$ is the restriction of $B_t^{\tilde{Y}}(s)$ on $Z$,
then we have
$$b_0|_Z={1\over \sqrt{\pi}}c_0.\eqno (5.36)$$

By (5.33), (5.36), one then gets
$$\int_Y a_0-\int_Y a_0^E={1\over \sqrt{\pi}}\int_Z c_0\ {\rm in}\ 
\Omega (B)/d\Omega (B).\eqno (5.37)$$

Now again using the results in Section 2, one has
$${1\over \sqrt{\pi}}\int_0^1(\int_Z c_0)ds=-\hat{\eta} (D^{Z,E},P_1,B^{Z,E})+
\hat{\eta}(D^E,P_0,B^Z)\ {\rm in}\ \Omega (B)/d\Omega (B),\eqno (5.38)$$
where $\{ P_s\}_{s\in [0,1]}$ is a (in fact any) curve of spectral sections
of $\{ D(s)=(1-s)D^Z+sD^{Z,E}\}_{ s\in [0,1]}$.

Combining (5.37), (5.38), we get
$$\int_Y a_0-\hat{\eta} (D^Z,P_0,B^Z)=\int_Ya_0^E-\hat{\eta}(D^{Z,E},P_1,B^{Z,E})
\ {\rm in}\ \Omega (B)/d\Omega (B).\eqno (5.39)$$

  From (5.39), one deduces that
$$\int_Y a_0-\hat{\eta} (D^Z,P,B^Z)=\int_Y a_0^E-\hat{\eta}(D^{Z,E},Q,B^{Z,E})$$
$$+\hat{\eta}(D^{Z,E},Q,B^{Z,E})-\hat{\eta}(D^{Z,E},P_1,B^{Z,E})
+\hat{\eta}(D^Z,P_0,B^Z)-\hat{\eta}(D^Z,P,B^Z)$$
$$={\rm ch}({\rm ind}\, (D^{Y,E},Q))+{\rm ch}([Q-P_1]+[P_0-P])$$
$$={\rm ch}({\rm ind}\, (D^{Y,E},Q))-{\rm ch}({\rm sf}\{ (D^{Z,E},Q),(D^Z,P)\} )$$
$$={\rm ch}({\rm ind}\, (D^Y,P)),\eqno (5.40)$$
where we have used results proved in this section and in Section 2.

The proof of Theorem 5.3 is completed. \qed

$\ $

{\bf Remark 5.4}. (5.25) may be viewed as an analogue in the family case of the
general index theorem of Atiyah-Patodi-Singer [APS2] for manifolds with
doundary.

X. Dai

Department of Mathematics, 
University of Southern  California, 
Los Angles, CA  90089 

{\em Current Address}: Math, UCSB, Santa Barbara, CA 93106

 xdai@math.usc.edu

W. Zhang

 Nankai Institute of Mathematics, Tianjin, 300071, P.R. China

\end{document}